\documentclass[
leqno,          
11pt,           
a4paper         
]{amsart}

\usepackage[colorlinks,citecolor=blue,urlcolor=blue]{hyperref}          
\usepackage{courier}                                                    
\usepackage{amssymb, amsmath, amsfonts, amsthm}                         
\usepackage{mathtools, cite, enumerate, rotating, framed, lipsum}       
\usepackage{fancybox, bbm, soul}                                        
\usepackage[dvipsnames]{xcolor}                                         
\usepackage[english]{babel}                                     
\usepackage[T1]{fontenc}                                                
\usepackage[cp1250]{inputenc}                                           

\numberwithin{equation}{section}                        

\newcommand{\thmcount}{thmc}                 
\newcounter{specialcounter}

\newtheorem{Thm}[\thmcount]{Theorem}
\newtheorem{Sthm}[specialcounter]{Theorem}
\newtheorem{Cor}[\thmcount]{Corollary}
\newtheorem{Lem}[\thmcount]{Lemma}
\newtheorem{Prop}[\thmcount]{Proposition}
\newtheorem{Rem}[\thmcount]{Remark}
\newtheorem{Defn}[\thmcount]{Definition}
\newtheorem{Ex}[\thmcount]{Example}
\newtheorem{Asu}[\thmcount]{Assumption}
\newtheorem{Sol}[\thmcount]{Solution}
\newtheorem*{Thmx}{Theorem}
\newtheorem*{Corx}{Corollary}
\newtheorem*{Lemx}{Lemma}
\newtheorem*{Propx}{Proposition}
\newtheorem*{Remx}{Remark}
\newtheorem*{Defnx}{Definition}
\newtheorem*{Exx}{Example}
\newtheorem*{Asux}{Assumption}
\newtheorem*{Solx}{Solution}

\newcommand \eq[1]{\begin{equation} #1 \end{equation}}
\newcommand \eqx[1]{\begin{equation*}  #1 \end{equation*}}
\newcommand \al[1]{\begin{align} #1 \end{align}}
\newcommand \alx[1]{\begin{align*}  #1 \end{align*}}
\renewcommand \sp[1]{\begin{equation} \begin{split} #1 \end{split} \end{equation}}
\newcommand \spx[1]{\begin{equation*} \begin{split} #1 \end{split} \end{equation*}}

\newcommand{\thm}[2]{\begin{Thm} \label{#1} #2 \end{Thm}}

\newcommand{\lem}[2]{\begin{Lem} \label{#1} #2 \end{Lem}}

\newcommand{\cor}[2]{\begin{Cor} \label{#1} #2 \end{Cor}}

\newcommand{\prop}[2]{\begin{Prop} \label{#1} #2 \end{Prop}}

\newcommand{\rem}[2]{\begin{Rem} \label{#1} #2 \end{Rem}}

\newcommand{\pr}[1]{\begin{proof} #1 \end{proof}}

\renewcommand{\a}{\alpha}                \newcommand{\e}{\varepsilon}
\newcommand{\w}{\omega}                 
        \newcommand{\la}{\lambda}

\newcommand{\RR}{\mathbb{R}}

\newcommand{\supp}{\mathrm{supp}}
\newcommand{\8}{\infty}
\renewcommand{\d}{\partial}
\renewcommand{\(}{\left(}
\renewcommand{\)}{\right)}

\newcommand{\Rd}{{\RR^d}}

\renewcommand{\rm}[1]{\mathrm{#1}}
\newcommand{\wt}[1]{\widetilde{#1}}

\newcommand{\sumk}{\sum_{k=1}^{\8}}

\newcommand{\abs}[1]{\left| #1 \right|}
\newcommand{\set}[1]{\left\{ #1 \right\}}
\newcommand{\norm}[1]{\left\| #1 \right\|}
\newcommand{\expr}[1]{\( #1 \)}
\newcommand{\eee}[1]{\( #1 \)}

\newcommand{\eeee}[1]{\left[ #1 \right]}
\newcommand{\lab}[1]{\label{#1}}

\newcounter{comcount}

\newcommand{\sleq}{\lesssim}
\newcommand{\sgeq}{\gtrsim}

\def\hadw{H^1_{at}(d,\w\mu)}
\def\hads{H^1_{at}(\rho,\sigma)}
\def\had1{H^1_{at}(d,\mu)}
\def\hadd{H^1_{at}(\wt{d}, \mu)}

\def\xy{(x,y)}

\title[  Hardy spaces for semigroups with Gaussian bounds ]{  Hardy spaces for semigroups with Gaussian bounds }

\author[ Jacek Dziuba\'nski ]{ Jacek Dziuba\'nski }
\address{
Jacek Dziuba\'nski \newline
\indent Instytut Matematyczny, Uniwersytet Wroc\l awski \newline
\indent pl. Grunwaldzki 2/4, 50-384 Wroc\l aw, Poland }
\email{jacek.dziubanski@math.uni.wroc.pl }

\author[ Marcin Preisner ]{ Marcin Preisner }
\address{
Marcin Preisner \newline
\indent Instytut Matematyczny, Uniwersytet Wroc\l awski \newline
\indent pl. Grunwaldzki 2/4, 50-384 Wroc\l aw, Poland }
\email{marcin.preisner@uwr.edu.pl }

\thanks{ The research  supported by  the Polish National Science Center (Narodowe Centrum Nauki) grant No. DEC--2012/05/B/ST1/00672}

\begin{document}

\begin{abstract}
Let $T_t=e^{-tL}$ be a semigroup of self-adjoint  linear operators acting on $L^2(X,\mu)$, where $(X,d,\mu)$ is a space of homogeneous type. We assume that $T_t$ has an integral kernel $T_t(x,y)$ which satisfies the upper and lower Gaussian bounds:
$$ \frac{C_1}{\mu(B(x,\sqrt{t}))} \exp\eee{{-c_1d(x,y)^2\slash t}}\leq T_t(x,y)\leq \frac{C_2}{\mu(B(x,\sqrt{t}))} \exp\eee{{-c_2 d(x,y)^2\slash t}}.$$
By definition, $f$ belongs to  $H^1_L$ if  $\| f\|_{H^1_L}=\|\sup_{t>0}|T_tf(x)|\|_{L^1(X,\mu)} <\infty$.
We prove that there is a function $\w(x)$, $0<c\le \w(x)\le C$, such that $H^1_L$ admits an atomic decomposition with atoms satisfying: ${\rm{ supp}}\, a\subset B$, $\| a\|_{L^\infty} \leq \mu(B)^{-1}$, and the weighted cancellation condition $\int a(x)\w(x)d\mu(x)=0$.
\end{abstract}

\maketitle

\section{Introduction}\lab{Intro}

Let $(X,d)$ be a metric space equipped with a non-negative Borel measure $\mu$. We shall assume that $\mu(X)=\8$ and $0<\mu(B(x,r))<\8$, $r>0$, where
    $
    B(x,r)= \set{y\in X \ | \ d(x,y)\leq r}
    $
 denotes the closed ball centered at $x$ and radius $r$. Suppose $(X,d,\mu)$ is a space of homogeneous type in the sense of Coifman-Weiss \cite{CoifmanWeiss_BullAMS}, which means that the doubling condition holds, namely: there is $C>0$ such that $\mu(B(x,2r))\leq C(\mu (B(x,r))$ for $r>0$ and $x\in X$. It is well-known that the doubling condition implies  that there are $q>0$ and $C_d>0$  such that
\begin{equation}\label{doubling}
\mu (B(x,sr))\leq C_d s^q\mu (B(x,r)) \ \ \text{for }  x\in X, \  r>0,  \  \text{and } s\geq 1.
\end{equation}

 Suppose that $L$ is a non-negative densely defined  self-adjoint linear operator on $L^2(X,\mu)$. Let $T_t=e^{-tL}$, $t>0$, denote the semigroup of linear operators generated by $-L$. We impose that there exists $T_t(x,y)$, such that
    \eq{\lab{int_rep}
    T_t f(x) = \int_X T_t(x,y) f(y) \, d\mu(y).
    }
Moreover, we assume the following lower and upper Gaussian bounds, that is, there are constants $c_1\geq c_2>0$ and $C_0>0$ such that
\begin{equation}\label{gauss12}
\frac{C_0^{-1}}{\mu(B(x,\sqrt{t}))} \exp\eee{-\frac{c_1d(x,y)^2}{t}} \leq T_t(x,y)\leq \frac{C_0}{\mu(B(x,\sqrt{t}))} \exp\eee{-\frac{c_2d(x,y)^2}{t}}
\end{equation}
for $t>0$ and a.e. $x,y\in X$. It is well-known that \eqref{gauss12} implies
    \eq{\label{gauss1212}
    \abs{\frac{\d^n}{\d t^n} T_t(x,y)}\leq \frac{C_n t^{-n}}{\mu(B(x,\sqrt{t}))} \exp\eee{-\frac{c'_n d(x,y)^2}{t}},
    }
for $t>0$ and a.e. $x,y\in X$  (for this fact see e.g. \cite[(7.1)]{Hofmann_Memoirs}, \cite{Davies}, \cite{Ouhabaz}).

The Hardy space $H^1(L)$ related to $L$ is defined by by means of the maximal function of the semigroup $T_t$, namely
    \eqx{
    H^1(L) := \set{f\in L^1(X,\mu) \ \Big| \ \norm{f}_{H^1(L)}:= \norm{\sup_{t>0}\abs{T_t f }}_{L^1(X,\mu)} <\8 }.
    }

On the other hand, we define atomic Hardy spaces as follows. Suppose we have a space $X$ with a doubling measure $\sigma$ and a quasi-metric $\rho$. We call a function $a$ an $(\rho,\sigma)$-atom if there exists a ball $B = B_\rho(x_0,r) := \set{x\in X \ | \ \rho(x,x_0)\leq r}$ such that:
$\supp \ a \subseteq B$, $\norm{a}_\8 \leq \sigma(B)^{-1}$, and
    \eqx{
    \int_B a(x)\, d\sigma(x) = 0.
    }
By definition, a function $f\in L^1(X, \sigma)$ belongs to $\hads$, if there exist $(\rho,\sigma)$-atoms $a_k$ and complex numbers $\la_k$, such that
\eqx{
f = \sum_{k=1}^\8 \la_k a_k \ \ \ \text{ and } \ \ \ \sumk \abs{\la_k} <\8.
}
If such sequences exist, we define the norm $\norm{f}_{\hads}$ to be the infimum of $\sumk \abs{\la_k}$ in the above presentations of $f$. Notice that in this paragraph we have changed the notation. This is because in the article we will use different metrics and measures.

For $x\in X$ let $\Phi_x:(0,\8) \to (0,\8)$ be a non-decreasing function defined by
\eq{\label{Psi}
\Phi_x(t) = \mu(B(x,t)).
}
The first main result of this paper is the following.
    \thm{mainthm1}{
    Suppose that $(X,d,\mu)$ satisfies \eqref{doubling} and assume that for each $x\in X$ the function $\Phi_x$ is a bijection on $(0,\8)$. Let an operator $L$ be given such that the semigroup $T_t$ satisfies \eqref{gauss12}. Then there exist a constant $C>0$ and a function $\w$ on X, $0< C^{-1} \leq \w(x) \leq C$, such that the spaces $H^1(L)$ and $\hadw$ coincide and the corresponding norms are equivalent,
        \eqx{
        C^{-1} \norm{f}_{\hadw} \leq \norm{f}_{H^1(L)} \leq C \norm{f}_{\hadw}.
        }
    Moreover, $\w$ is $L$-harmonic, that is $T_t \w=\w$ for $t>0$.
    }

Let us notice that the assumption on $\Phi_x$ implies that $\mu(X) = \8$ and $\mu$ is non-atomic. This will be used later on. By definition, we call a semigroup {\it conservative} if
    \eq{\label{conservative}
    \int_X T_t(x,y) \, d\mu(y)=1
    }
for $t>0$ and $x\in X$.

\thm{mainthm2}{
Suppose that $(X,d,\mu)$ satisfies \eqref{doubling} and assume that for each $x\in X$ the function $\Phi_x$ is a bijection on $(0,\8)$. Let an operator $L$ be given such that the semigroup $T_t$ satisfies \eqref{gauss12} and \eqref{conservative}. Then the spaces $H^1(L)$ and $\had1$ coincide and the corresponding norms are equivalent, i.e. there exists $C>0$ such that
    \eqx{
    C^{-1} \norm{f}_{\had1} \leq \norm{f}_{H^1(L)} \leq C \norm{f}_{\had1},
    }
}

It appears that Theorem \ref{mainthm2} is equivalent to Theorem \ref{mainthm1}, see Section \ref{Doob}. Let us also emphasise that we do not require any regularity  conditions on the kernels $T_t(x,y)$. However, it turns out that \eqref{gauss12} implies crucial for this paper H\"older type estimates on $T_t(x,y)$. This will be discussed in Section  4 (see Theorem \ref{Lip_Gauss} and Corollary \ref{Holder_Gauss}).

The theory of the classical Hardy spaces on the Euclidean spaces $\mathbb R^n$ has its origin in studying holomorphic functions of one variable in the upper half-plane. The reader is referred to the original works of Stein and Weiss \cite{Stein_Weiss}, Burkholder et al. \cite{Burkholder_Gundy_Silverstein}, Fefferman and Stein \cite{Fefferman_Stein}. Very important contribution to the theory is atomic decompositions of the elements of the $H^p$ spaces  proved by  Coifman \cite{Coifman_Studia} in the one dimensional case and then by Latter \cite{Latter_Studia} for $H^p(\mathbb R^n)$ . The theory was then extended to the space of homogeneous type (see e.g., \cite{CoifmanWeiss_BullAMS}, \cite{Macias_Segovia_Advances}, \cite{Uchiyama}). For more information concerning  the classical Hardy spaces, their characterisations and historical comments we refer the reader to Stein \cite{Stein}. A very general approach to the theory of Hardy spaces associated with semigroups of linear operators satisfying the Davies-Gaffney estimates was introduced by Hofmann et al. \cite{Hofmann_Memoirs} (see also \cite{Auscher_unpublished}, \cite{Song_Yan_homogeneous}). Let us point out that the classical Hardy spaces can be thought as those associated with the classical heat semigroup $e^{t\Delta}$. Finally we want to remark that the present paper takes motivation from \cite{DZ_JFAA} and \cite{DZ_Potential_2014}, where the authors studied $H^1$ spaces associated with Schr\"odinger operators $-\Delta +V$ on $\mathbb R^n$, $n\geq 3$, with Green bounded  potentials $V\geq 0$.

In what follows $C$ and $c$ denotes different constants that may depend on $C_d, q, C_0, c_1, c_2$. By $U \sleq V$ we understand $U\leq CV$ and $U \simeq V$ means $C^{-1}V \leq U \leq CV$.

The paper is organized as follows. In Section \ref{harmonic} we prove that the estimates \eqref{gauss12} imply existence of $L$-harmonic function $\w$ such that $0<C^{-1}\leq \w(x)\leq C$. The equivalence of Theorems \ref{mainthm1} and \ref{mainthm2} is a consequence of the Doob transform, see Section \ref{Doob}. Then a proof of Theorem \ref{mainthm2} is given in a few steps. First, in  Section \ref{sec_Lip} we prove H\"older-type estimates for a conservative semigroup. Then, in Section \ref{sec_new_met} we introduce a new quasi-metric $\wt{d}$ and study its properties. In section \ref{Uchiyama} we apply a theorem of Uchiyama on the space $(X,\wt{d}, \mu)$ to  complete the proof of Theorem \ref{mainthm2}. Finally, in Section \ref{sec_ex} we provide some examples of semigroups that satisfy assumptions of Theorem \ref{mainthm1}.

\section{Gaussian estimates and bounded harmonic functions}\label{harmonic}

In this section we assume that the semigroup $T_t$ satisfies \eqref{gauss12}. Clearly,
$$ C^{-1} \leq \int_X T_t(x,y)d\mu(y)\leq C$$
for all $t>0$ and a.e. $x\in X$. For a positive integer $n$ define
$$ \w_n(x)=\frac{1}{n}\int_0^n\int_X T_s(x,y)d\mu (y)\, ds.$$
Then,
$$ C^{-1}\leq \w_n(x) \leq C.$$

Recall that a metric space with the doubling condition is  separable, then so is $L^1(X,\mu)$.  Using the Banach-Alaoglu theorem for $L^\8(X,\mu)$ there exists a subsequence $n_k$ and $\w\in L^\infty(X,\mu)$, such that $\w_{n_k} \to \w$ in $*$-weak topology. Obviously, 
$$ C^{-1} \leq \w(y)\leq C.$$

Our goal is to prove that $T_t \w (x) = \w(x)$ for $t>0$. To this end, we write
\begin{equation*}
\begin{split}
\int_X T_t(x,y)\w(y)\,d\mu(y)&=\lim_{k\to\infty} \int_X T_t(x,y) \w_{n_k}(y)\, d\mu (y)\\
&=\lim_{k\to\infty} \frac{1}{n_k}  \int_0^{n_k} \int_X T_{t+s}(x,z) \, d\mu(z)\, ds \\
&=\lim_{k\to\infty} \w_{n_k}(x)
 + \lim_{k\to\infty} \frac{1}{n_k}  \int_{n_k}^{n_k+t} \int_X T_{s}(x,z) \, d\mu(z)\, ds \\
&\ \ \ -\lim_{k\to\infty} \frac{1}{n_k}  \int_{0}^{t} \int_X T_{s}(x,z) \, d\mu(z)\, ds. \\
\end{split}
\end{equation*}
Since the last two limits tend to zero, as $k\to \8$, we obtain
\begin{equation}\label{eq1}  \int_X T_t(x,y)\w(y)\,d\mu(y)=
\lim_{k\to\infty} \w_{n_k}(x).\end{equation}
From \eqref{eq1} we get that $\lim_{k\to\infty} \w_{n_k}(x)$ exists for a.e. $x\in X$ and the limit has to be $\w(x)$. Moreover, $T_t \w(x) = \w(x)$. Thus we have proved the following proposition.

\prop{coro_harm}{
Assume that a semigroup $T_t$ satisfies \eqref{gauss12}. Then there exists a function $\w$ and $C>0$ such that $0<C^{-1}\leq \w(x)\leq C$ and $T_t \w(x) = \w(x)$ for every $t>0$.
}

\section{Doob transform}\label{Doob}

In this section we work in a slightly more general scheme. Let $(X,\mu)$ be a $\sigma$-finite measure space and $L$ be a self-adjoint operator on $L^2(X,\mu)$. We assume that the strongly continuous semigroup $T_t=\exp(-tL)$ admits a non-negative integral kernel $T_t(x,y)$, so that
    \eqx{
    T_t f(x) = \int_X T_t(x,y) f(y) \, d\mu(y).
    }
Moreover, assume that there exists $\w$ satisfying $0<C^{-1} \leq \w(x) \leq C$, that is $L$-harmonic. Namely, for every $t>0$ one has
$$\int_X T_t(x,y)\w(y)\, d\mu (y)=\w(x)$$
for a.e. $x\in X$. Obviously, this implies that $\sup_{x\in X, t>0}\int T_t(x,y)\, d\mu (y)\leq C$.

Define a new measure $d\nu (x)=\w^2(x)\, d\mu(x)$ and a new kernel
\eq{\label{TK}
K_t(x,y)=\frac{T_t(x,y)}{\w(x)\w(y)}.
}
The semigroup $K_t$ given by
    \eqx{
    K_t f(x) = \int_X K_t(x,y) f(y) \, d\nu(y)
    }
is a strongly continuous semigroup of self-adjoint integral operators on $L^2(X,\nu)$. The mapping $L^2(X, \mu)\ni f\mapsto \w^{-1} f \in L^2(X,\nu)$ is an isometric isomorphism with the inverse $L^2(X, \nu)\ni g\mapsto  \w g\in L^2(X,\mu)$. Clearly,
$$ K_tg(x)=\w(x)^{-1} T_t(\w g)(x).$$
Moreover, the positive integral kernel $K_t(x,y)$ is conservative, that is,
$$  \int_X K_t(x,y)\, d\nu (y)= \w(x)^{-1} T_t\w(x) =1.$$
Thus the above change of measure and operators, which is called Doob's transform (see e.g., \cite{Grigor'yan_HeatKernels}), conjugates the semigroup $T_t$ with the conservative semigroup $K_t$.

It is clear that the operators $K_t$ are contractions on $L^1(X,\nu)$. Consequently, $K_t$ is a strongly continuous semigroup of linear operators  on $L^1(X,\nu)$. To see this, it suffices to show that $\lim_{t\to 0} \| K_t\chi_A-\chi_A\|_{L^1(X,\nu)}=0$ for any measurable set $A$ of finite measure. We have that
$$ \lim_{t\to 0}\int_A | K_t\chi_A -\chi_A|d\nu \leq \lim_{t\to 0} \| K_t\chi_A -\chi_A\|_{L^2(X,\nu)}\nu (A)^{1\slash 2}=0.$$
On the other hand $\int_X K_t\chi_A(x)\, d\nu(x)= \nu(A)=\| \chi_A\|_{L^1(X,\nu)}$. Hence,
\begin{equation*}\begin{split}
\int_{A^c} |K_t\chi_A(x)|d\nu(x) & = \int_{A^c} K_t\chi_A(x)d\nu(x)=\int_X K_t\chi_A(x)d\nu (x)-\int_{A} K_t\chi_A(x)\, d\nu (x)\\
&=\int_A \chi_A(x)d\nu (x)-\int_{A} K_t\chi_A(x)\, d\nu (x) \to 0 \ \ {\rm as} \ t\to 0,\\
\end{split}\end{equation*}
which completes the proof of the strong continuity of $K_t$ on $L^1(X,\nu)$.

Further, we easily see that the semigroup  $T_t$ is strongly continuous on $L^1(X,\mu)$. Indeed, if $f\in L^1(X,\mu)$ then $g=\w^{-1} f\in L^1(X,\nu)$ and
\begin{equation*}
\begin{split}
\| T_tf -f\|_{L^1(X,\mu)} &=\int_X\abs{K_t g(x)-g(x)} \w(x)^{-1}  \,d\nu(x)\\
    &\sleq \int_X\abs{K_t g(x)-g(x)}\,d\nu(x) \to 0 \ \ {\rm as} \ t\to 0.
\end{split}
\end{equation*}

Now we discuss the equivalence of Theorems \ref{mainthm1} and \ref{mainthm2}. Let $T_t$ and $K_t$ be the semigroups related by \eqref{TK} with generators $L$ and $R$, respectively. It easily follows from the Doob transform, that
$f\in H^1(L)$ if and only if
$\w^{-1}f\in H^1(R)$
and $\| f\|_{H^1(L)}\simeq \| \w^{-1}f\|_{H^1(R)}$. In other words
$$ {H^1(L)}\ni f\mapsto \w^{-1}f\in H^1(R)$$
is an isomorphism of the spaces.

Assume that the space $H^1(R)$ admits an atomic decomposition with atoms that satisfy the  cancellation condition with respect to the measure $\nu$, that is every $g\in H^1(R)$ can be written as $g=\sum\lambda_j a_j$ with $\sum_{j} |\lambda_j|\simeq \| g\|_{H^1(R)}$ and $a_j$ are atoms with the property $\int a_j\, d\nu =0$. Then every $f\in {H^1(L)}$ admits atomic decomposition $f=\sum \lambda_j b_j$ with atoms $b_j$ that satisfy cancellation condition $\int b_j\,  \w\, d\mu=0$.

\section{H\"older-type estimates on the semigroups}\label{sec_Lip}

In this section we consider a conservative  semigroup $T_t$ having an integral kernel $T_t(x,y)$ that satisfies the upper and lower Gaussian bounds \eqref{gauss12}.  Let
    \begin{equation}\label{subord}
    P_t(x,y)= \pi^{-1/2}\int_0^{\infty} e^{-u}T_{t^2\slash (4u)}(x,y) \frac{du}{\sqrt{u}}
    \end{equation}
be the kernels of the subordinate semigroup $P_t = e^{-t\sqrt{L}}$.

\thm{Lip_Poiss}{
Assume that the semigroup $T_t$ satisfies \eqref{gauss12} and \eqref{conservative}. Then there is a constant $\alpha>0$ such that
    \begin{equation}
    |P_t(x,y)-P_t(x,z)|\sleq \Big(\frac{d(y,z)}{t}\Big)^\alpha P_t(x,y)
    \end{equation}
    whenever $d(y,z)\leq t$.
}

\thm{Lip_Gauss}{
Assume that the semigroup $T_t$ satisfies \eqref{gauss12} and \eqref{conservative}. Then there are constants $\alpha,c>0$ such that
    \begin{equation}\label{gauss13}
    |T_t(x,y)-T_t(x,z)|\sleq \mu(B(x,\sqrt{t}))^{-1} \eee{\frac{d(y,z)}{\sqrt{t}}}^\alpha \exp\eee{-\frac{cd(x,y)^2}{ t}}
    \end{equation}
whenever $d(y,z)\leq \sqrt{t}$.
}

H\"older regularity of semigroups satisfying Gaussian bounds was considering in various settings by many authors.
We refer the reader to Grigor'yan \cite{Grigor'yan_HeatKernels}, Hebisch and Saloff-Coste \cite{Hebisch_Saloff-Coste_Grenoble}, Saloff-Coste \cite{Saloff-Coste}, Gyrya and Saloff-Coste \cite{Gyrya_Saloff-Coste}, Bernicot, Coulhon, and Frey \cite{Bernicot_Coulhon_Frey}, and references therein.
Here we present a short alternative  proof of \eqref{gauss13}. To this end we  shall first prove some auxiliary propositions and then  Theorem \ref{Lip_Poiss}. Finally, at the end of this section, we shall make use of functional calculi deducing Theorem \ref{Lip_Gauss} from Theorem \ref{Lip_Poiss}.

\prop{Poiss_new}{
For $x,y \in X$ and $t>0$ we have
    \eq{\label{poiss_new}
    P_t(x,y) \simeq \frac{1}{\mu\eee{B\(x, t+d(x,y)\)}} \ \frac{t}{t+d(x,y)}.
    }
}

\pr{
Consider first the case  $d(x,y) \leq t$. Then $d(x,y)+t\simeq t$. The upper bound follows by \eqref{gauss12} and \eqref{doubling}. Indeed,
        \begin{equation*}\begin{split}
        P_t(x,y) & \sleq \int_{0}^{1 \slash 4} \mu\Big(B\Big( x,\frac{t}{2\sqrt{u}}\Big)\Big)^{-1} \frac{du}{\sqrt{u}} + \mu(B(x,t))^{-1} \int_{1\slash 4}^\infty  \frac{e^{-u}\mu(B(x,t))}{\mu\(B\( x,\frac{t}{2\sqrt{u}}\)\)} \frac{du}{\sqrt{u}}\\
        &\sleq \mu(B( x,t))^{-1}\int_{0}^{1 \slash 4} \frac{du}{\sqrt{u}} + \mu(B(x,t))^{-1} \int_{1\slash 4}^\infty  e^{-u} (2\sqrt{u})^q \frac{du}{\sqrt{u}}\\
        &\sleq \mu(B( x,t))^{-1}.
        \end{split}\end{equation*}
        Also, \eqref{gauss12} implies lower bounds, since
        \begin{equation*}
        P_t(x,y)\sgeq \int_{1\slash 4}^1 \frac{e^{-u} \exp\eee {- 4c_1 u}}{\mu\(B\( x,\frac{t}{2\sqrt{u}}\)\)} \frac{du}{\sqrt{u}}\sgeq \mu(B(x,t))^{-1}.
        \end{equation*}

Let us now turn to the case $d(x,y)\geq t$. Then $d(x,y)+t\simeq d(x,y)$. Using \eqref{gauss12}, we have
    \begin{equation}\label{int_7}\begin{split}
    P_t(x,y) &\sleq \int_0^{\8} \frac{e^{-u} \exp\eee{-4c_2 ud(x,y)^2\slash t^2}}{\mu \(B\(x,\frac{t}{2\sqrt{u}}\)\)}\frac{du}{\sqrt{u}} = \int_0^{t^2\slash d(x,y)^2} +
    \int_{t^2\slash d(x,y)^2}^\8 = (J_1) + (J_2).
    \end{split}\end{equation}
Moreover, by \eqref{doubling},
    \begin{equation}\label{int_1}\begin{split}
    (J_1) &\simeq \int_0^{t^2\slash d(x,y)^2} \frac{1}{\mu \(B\(x,\frac{t}{2\sqrt{u}}\)\)}\frac{du}{\sqrt{u}}\\
    &\sleq \int_0^{t^2\slash d(x,y)^2} \frac{1}{\mu \(B\(x,\frac{d(x,y)}{2}\)\)}\frac{du}{\sqrt{u}}
    \simeq \frac{1}{\mu\(B(x,d(x,y))\)} \ \frac{ t}{d(x,y)}
    \end{split}\end{equation}
and
    \begin{equation}\label{int_3}\begin{split}
    (J_2) &\sleq \frac{1}{\mu\(B(x,d(x,y))\)} \int_{t^2\slash d(x,y)^2}^{\infty}  \exp\eee{-4c_2 ud(x,y)^2\slash t^2}
    \Big(\frac{2d(x,y)\sqrt{u}}{t}\Big)^q \frac{du}{\sqrt{u}}\\
    &\simeq \frac{1}{\mu\(B(x,d(x,y))\)} \ \frac{t}{d(x,y)}. \end{split}\end{equation}

The estimates \eqref{int_7}--\eqref{int_3} give upper estimate. For the lower estimate, recall that $d(x,y)\geq t$ and observe that
    \begin{equation*}\begin{split}
    P_t(x,y) &\sgeq \int_{t^2\slash d(x,y)^2}^{2t^2\slash d(x,y)^2}
     \frac{ e^{-u} \exp\eee{-4c_1ud(x,y)^2\slash t^2}}{\mu \(B\(x,\frac{t}{2\sqrt{u}}\)\)}\frac{du}{\sqrt{u}}\\
     &\simeq \frac{1}{{\mu(B(x,d(x,y)\slash 2))}} \int_{t^2/d(x,y)^2}^{2t^2/d(x,y)^2}\frac{du}{\sqrt{u}}\\
     &\simeq \frac{1}{\mu (B(x,d(x,y)))} \ \frac{t}{d(x,y)}.
    \end{split}\end{equation*}
}

\cor{coro_2}{
     There is a constant $C>0$ such that if $d(y,z)\leq t$, then
        \begin{equation}\label{between}
        C^{-1} \leq \frac{P_t(x,y)}{P_t(x,z)}\leq C.
        \end{equation}
}

\pr{
        The corollary is a simple consequence of \eqref{poiss_new}. For the proof one may consider two cases: $d(x,y)\leq 2t$ (then $d(x,z)\leq  3t$) and $d(x,y)>2 t$ (then $d(x,y)\simeq d(x,z)$).
}

\prop{shrink}{
There exists a constant $\gamma\in(0,1)$ such that the following statement holds:\\ \noindent
If there are $y,z,t_0,a_1,b_1>0$ given such that $d(y,z)<t_0$ and for all $x\in X$ we have
    \begin{equation}\label{base}
    a_1\leq \frac{P_{t_0}(x,y)}{P_{t_0}(x,z)} \leq b_1,
    \end{equation}
then there is a subinterval $[a_2,b_2] \subseteq [a_1,b_1]$ such that $b_2-a_2 = \gamma(b_1-a_1)$ and for all $t\geq 2t_0$ and all $x\in X$ one has
    $$
    a_2\leq \frac{P_{t}(x,y)}{P_{t}(x,z)} \leq b_2.
    $$
}

\pr{
The proof, which takes some ideas from \cite{Jerison_Kenig}, is an adapted version of the proof of \cite[Prop. 3.1]{DZ_JFAA}. For the reader convenience we present the details. Let $m=(a_1+b_1)/2$ and $\theta = (b_1-a_1)/(a_1+b_1) \in(0,1)$, so that $a_1=(1-\theta)m$ and $b_1=(1+\theta)m$. Define
    $$
    \Omega_+ = \set{x\in X \ \Big| \ m \leq \frac{P_{t_0}(x,y) }{P_{t_0}(x,z)}\leq b_1}, \qquad \Omega_- = X \setminus \Omega_+.
    $$
Obviously, either $\mu(\Omega_- \cap B(z,t_0)) \geq \mu(B(z,t_0))/2$ or $\mu(\Omega_+ \cap B(z,t_0)) \geq \mu(B(z,t_0))/2$. Here we shall assume that the latter holds. The proof in the opposite case is similar.
Denote $B = B(z,t_0)$, $s=t-t_0 \geq t_0$. For $x\in X$, we have
    \begin{align*}
    P_t (x,y) &\geq m(1-\theta) \int_{\Omega_-} P_s(x,w) P_{t_0}(w,z)\, d\mu(w) + m \int_{\Omega_+} P_s(x,w) P_{t_0} (w,z)\, d\mu(w) \\
    &= m(1-\theta) P_t(x,z) + m\theta \int_{\Omega_+} P_s(x,w) P_{t_0} (w,z)\, d\mu(w)\\
    &\geq m(1-\theta) P_t(x,z) + m\theta \int_{\Omega_+\cap B} P_s(x,w) P_{t_0} (w,z)\, d\mu(w)\\
    &\geq m(1-\theta) P_t(x,z) + m\theta \frac{\mu(B)}{2} \inf_{w\in B} P_s (x,w) \inf_{w\in B} P_{t_0} (w,z) = (J).
    \end{align*}
Notice that $d(w,z)\leq t_0 \leq s$ and, by Corollary \ref{coro_2},
    $$
    \inf_{w\in B} P_s (x,w) \simeq P_s(x,z).
    $$
Since $t\simeq s$, Proposition \ref{Poiss_new} implies
    $$
    P_s(x,z)\simeq P_t(x,z)
    $$
and
    $$
    \inf_{w\in B} P_{t_0} (w,z) \simeq \mu(B)^{-1}.
    $$
Therefore
    \begin{equation}\label{J}
    (J) \geq m\left((1-\theta) + \kappa \theta \right) P_t(x,z) = m\left(1-\theta(1-\kappa)\right)P_t(x,z),
    \end{equation}
where $\kappa\in(0,2)$. Moreover, from \eqref{base} and the semigroup property we easily get
    \begin{equation}\label{J2}
    P_t(x,y) \leq b_1 \int P_s(x,w) P_{t_0}(w,z) \, d\mu(w) \leq b_1 P_t(x,z).
    \end{equation}
Defining $\gamma = 1-\kappa/2 \in(0,1)$, $b_2=b_1$ and $a_2 = m\left((1-\theta) + \kappa \theta \right)$, we have
    $$
    b_2-a_2 =2m\theta\left(1-\kappa/2\right) = (b_1-a_1)(1-\kappa/2) = \gamma(b_1-a_1).
    $$
Now \eqref{J} together with \eqref{J2} give
    $$
    a_2 \leq \frac{P_t(x,y)}{P_t(x,z)} \leq b_2 \ \ \text{ for } x\in X.
    $$
}

\pr{[Proof of Thorem \ref{Lip_Poiss}]
Having Corollary \ref{coro_2} and Proposition \ref{shrink} proved, we follow arguments of \cite{DZ_JFAA} to obtain the theorem. By Corollary \ref{coro_2} there are $b_1>a_1>0$ such that for $y,z \in X$ and $t>0$ satisfying $d(y,z)<t$ we have
    $$
    a_1 \leq \frac{P_t(x,y)}{P_t(x,z)} \leq b_1
    $$
for all $x\in X$. From Proposition \ref{shrink} we deduce that there exists $\w(y,z)$ such that
    \begin{equation}\label{conv1}
    \lim_{t\to \infty} \frac{P_t(x,y)}{P_t(x,z)} = \w(y,z) \quad \textrm{uniformly in  }x\in X.
    \end{equation}
It follows from \eqref{subord} that $\int_X P_t(x,y)\, d\mu(y) = 1$. Recall that $P_t(x,y) = P_t(y,x)$. Using \eqref{conv1},
    $$
    1 = \int P_t(y,x) \, d\mu(x)
    = \int \frac{P_t(y,x)}{P_t(z,x)} {P_t(z,x)} \, d\mu(x) \xrightarrow[t\to \infty]{} \w(y,z).
    $$
Thus $\w(y,z) = 1$.

Assume that $d(y,z)<t$. Let $n\in \mathbb{N}$ be such that $d(y,z) \leq t2^{-n} < 2d(y,z)$. Set $t_0 = t2^{-n}$. Clearly, $d(y,z)\leq t_0$ and
    $$
    a_1 \leq \frac{P_{t_0}(x,y)}{P_{t_0}(x,z)} \leq b_1.
    $$
Observe that $n\simeq \log(t/d(y,z))$. Applying Proposition \ref{shrink} $n$-times we arrive at
    $$
    \left|\frac{P_t(x,y)}{P_t(x,z)}-1\right|\sleq \gamma^{n} \sleq \gamma^{c\log(t/d(y,z))} \sleq \left(\frac{d(y,z)}{t}\right)^\alpha.
    $$
with $\alpha>0$ and the proof of Theorem \ref{Lip_Poiss} is finished.
}

Finally, we devote the remaining part of this section for  deducing  Theorem \ref{Lip_Gauss} from Theorem \ref{Lip_Poiss}. This is done by using a functional calculi. First, we need some preparatory facts.
Recall that $q$ is a fixed constant satisfying \eqref{doubling}. By $W^{2,\sigma}(\RR)$ we denote the Sobolev space with the norm
    \eqx{
    \norm{f}_{W^{2,\sigma}(\RR)} = \eee{\int_\RR \eee{1+|\xi|^2}^\sigma \abs{\widehat{f}(\xi)}^2\, d\xi}^{1/2}.
    }

Let $\int_0^\infty \xi  dE_{\sqrt{L}}(\xi)$ be the  spectral resolution for $\sqrt{L}$. For a bounded function $m$ on $[0,\infty)$ the formula
$$ m(\sqrt{L})=\int_0^\infty m(\xi )\, dE_{\sqrt{L}}(\xi)$$
defines a bounded linear operator on $L^2(X,\mu)$.

Further we shall use the following lemma, whose proof based on finite speed propagation of the wave equation (see \cite{Coulhon_Sikora}) can be found in \cite{DP_Argentina}.

        \lem{LemaRevista}{\cite[Lemma 4.8]{DP_Argentina}
           Let $\kappa >1/ 2$, $\beta >0$. Then there exists a constant $C>0$ such that
         for every even function  $m\in W^{2,\beta\slash 2+\kappa}(\mathbb R)$ and every  $g\in
         L^2(X,\mu)$, $\mathrm{supp}\, g\subset B(y_0, r)$,
         we have
          $$ \int_{d(x,y_0)>2r} \abs{m(2^{-j}\sqrt{
          L})g(x)}^2\left(\frac{d(x,y_0)}{r}\right)^\beta d\mu(x)\leq C
          (r2^j)^{-\beta}
          \| m\|_{W^{2,\beta\slash 2+\kappa}(\RR)}^2\| g\|_{L^2(X,\mu)}^2$$
          for $j\in\mathbb Z$.
        }

\prop{prop112}{
Let $\beta>q$ and $\kappa >1\slash 2$. There is a constant $C>0$ such that for every $F\in W^{2, \beta \slash 2+\kappa}(\RR)$ with $\mathrm{supp}\, F\subset (1/2, 2)$ the integral kernels $F(2^{-j}\sqrt{L})(x,y)$ of the operators $F(2^{-j}\sqrt{L})$ satisfy
    \eqx{\label{ineq50}
    \int_X |F(2^{-j}\sqrt{L})(x,y)|\, d\mu (x)\leq C \| F\|_{W^{2, \beta/2+\kappa}(\RR)}  \ \ \text{ for } \ j\in\mathbb Z.
    }
}

\begin{proof}
For $y\in X$ and $j\in\mathbb Z$ set $U_0=B(y,2^{-j})$, $U_k=B(y,2^{k-j})\setminus B(y,2^{k-j-1})$, $k=1,2,...$ . Define $g_{j,k,y}(x)=T_{2^{-2j}}(x,y)\chi_{U_k}(x)$, $k=0,1,2,...$\,. Then, using \eqref{doubling} and \eqref{gauss12}, we have
\begin{equation}\label{ineq51}
\| g_{j,y,k}\|_{L^2(X,\mu)}\leq a_k \mu(U_0)^{-1\slash 2},
\end{equation}
where $a_k=C_0 \sqrt{C_{q}2^{kq}}\exp(-c_22^{2k-2})$ is a rapidly decreasing sequence.

Let $m(\xi)$ be the even extension of $e^{\xi^2}F(\xi)$. Obviously, $\norm{m}_{W^{2, \beta/2+\kappa}(\RR)} \simeq \norm{F}_{W^{2, \beta/2+\kappa}(\RR)}$.  Then,
$F(2^{-j}\sqrt{L})=m(2^{-j}\sqrt{L})T_{2^{-2j}}$,
and, consequently,
\begin{equation}\label{ineq51a}
F(2^{-j}\sqrt{L})(x,y)=\sum_{k= 0}^\infty m(2^{-j}\sqrt{L})g_{j,y,k}(x).
\end{equation}
By  the Cauchy-Schwartz inequality, \eqref{doubling}, and \eqref{ineq51} we get
\begin{equation}\label{ineq52}\begin{split}
  &\|m(2^{-j}\sqrt{L})g_{j,y,k}\|_{L^1(B(y,2^{k-j+1}),\mu)}
  \leq \mu(B(y,2^{k-j+1}))^{1\slash 2} \|m(2^{-j}\sqrt{L})g_{j,y,k}\|_{L^2(X,\mu)}\\
  &\sleq  \(\frac{\mu(B(y,2^{k-j+1}))}{\mu(B(y,2^{-j}))}\)^{1\slash 2} \mu(B(y,2^{-j}))^{1\slash 2} \|m\|_{L^\infty(\RR)}\|g_{j,y,k}\|_{L^2(X,\mu)}\\
  &\sleq 2^{kq} a_k \| m\|_{W^{2, \beta/2+\kappa}(\RR)}.
\end{split}\end{equation}
We turn to estimate $\|m(2^{-j}\sqrt{L})g_{j,y,k}\|_{L^1(B(y,2^{k-j+1})^c,\mu)}$. Utilizing  the Cauchy-Schwarz inequality and Lemma \ref{LemaRevista}, we obtain,
\begin{equation}\begin{split}\label{ineq53}
                           & \int_{d(x,y)>2^{k-j+1}}\abs{m(2^{-j}\sqrt{L})g_{j,y,k}(x)}\, d\mu(x)\\
                            &  \leq \eeee{ \int_{d(x,y)>2^{k-j+1}}\abs{m(2^{-j}\sqrt{L})g_{j,y,k}(x)}^2\(\frac{d(x,y)}{2^{k-j}}\)^\beta\, d\mu(x)}^{1\slash 2} \\
                            & \ \ \times
                            \eeee{ \int_{d(x,y)>2^{k-j+1}}\eee{\frac{d(x,y)}{2^{k-j}}}^{-\beta}\, d\mu(x)}^{1\slash 2} \\
                              & \sleq 2^{-k\beta\slash 2}  \| m\|_{W^{2,\beta\slash 2+\kappa}(\RR)} \|g_{j,y,k}\|_{L^2(X,\mu)} \eeee{ \int_{d(x,y)>2^{k-j+1}}\Big(\frac{d(x,y)}{2^{k-j}}\Big)^{-\beta}\, d\mu(x)}^{1\slash 2}.\\
                         \end{split}
\end{equation}
Recall that $\beta >q$, hence is not difficult to check that \eqref{doubling} leads to
\begin{equation}\label{ineq54}
   \int_{d(x,y)>2^{k-j+1}}\eee{\frac{d(x,y)}{2^{k-j}}}^{-\beta}\, d\mu(x)\sleq 2^{kq}\mu(U_0).
\end{equation}
Thus, by \eqref{ineq51}, \eqref{ineq53} and \eqref{ineq54}, we get
\begin{equation*}\begin{split}\label{ineq55}
                            \int_{d(x,y)>2^{k-j+1}}\abs{m(2^{-j}\sqrt{L})g_{j,y,k}(x)}\, d\mu(x)
                               \sleq 2^{-k\beta\slash 2}  2^{kq\slash 2}\| m\|_{W^{2,\beta\slash 2+\kappa}(\RR)} a_k,\\
                         \end{split}
\end{equation*}
which, combined with \eqref{ineq52} and \eqref{ineq51a}, completes the proof of the proposition.
\end{proof}

{For $t>0$ set
    $$
    \psi_t(\xi) = \exp\eee{\sqrt{t}\xi -t\xi^2}.
    $$

\lem{spec_mult}{
    The operators $\psi_t(\sqrt{L})$ have integral kernels $\Psi_t(x,y)$ that satisfy
        \begin{align}
        &\psi_t(\sqrt{L})f(x) = \int \Psi_t(x,y) f(y) \,d\mu(y), \notag \cr
        &\sup_{x\in X,\, t>0} \ \int |\Psi_t(x,y)|\, d\mu(y) \leq C. \label{kern_est}
        \end{align}
}

\pr{
    Denoting $\theta_t(\xi) = \exp\eee{-t\xi^2}(\exp\eee{\sqrt{t}\xi}-1)$, where $t,\xi>0$, we have
        $$
        \psi_t(\sqrt{L}) = T_t + \theta_t(\sqrt{L}).
        $$
    Clearly, by \eqref{gauss12}, $\sup_{t,x>0} \int T_t (x,y) \, d\mu(y) \leq C$. Thus we concentrate our attention on $\theta_t(\sqrt{L})$.

        Let $\eta \in C_c^\infty (1/2,2)$ be a partition of unity such that
        $$
        \theta_t(\xi) = \sum_{j\in {\mathbb{Z}}} \theta_t(\xi) \eta(2^{-j}\xi) = \sum_{j\in {\mathbb{Z}}} \theta_{t,j}(\xi).
        $$
    Denote $\wt{\theta_{t,j}}(\xi) = \theta_{t,j}(2^j \xi) = \eta(\xi) \theta_1(\sqrt{t} 2^j\xi)$. One can easily verify that for $n\in \mathbb{N}\cup\{0\}$ there are constants $C_n, c_n>0$ such that
        $$
        \Big\|\frac{d^n}{d\xi^n} \wt{\theta}_{t,j}\Big\|_{L^\infty(\RR)} \leq C_n
        \begin{cases}
        \sqrt{t}2^{j}\quad \ \,  &\text{for } j\leq -\log_2 \sqrt{t}\\
        \exp\eee{-c_n t2^{2j}}\qquad \  &\text{for } j\geq -\log_2 \sqrt{t}
        \end{cases}.
        $$
    In other words, for arbitrary $N$,
        $$
        \sum_{j\in {\mathbb{Z}}} \|\wt{\theta}_{t,j}\|_{W^{2,N}(\RR)}\leq C(N).
        $$

Using Proposition \ref{prop112} (with a fixed $N> q\slash 2+1$)  we get that the integral kernels $\Theta_{t,j}(x,y)$ of the operators $\theta_{t,j} (\sqrt{L}) = \wt{\theta}_{t,j}(2^{-j} \sqrt{L})$ satisfy
        $$
        \int |\Theta_{t,j}(x,y)| d\mu(y) \sleq  \|\wt{\theta}_{t,j}\|_{W^{2,N}(\RR)}.
        $$
    Therefore $\Theta_t (x,y) = \sum_{j\in {\mathbb{Z}}} \Theta_{t,j}(x,y)$ is the integral kernel of $\theta_t(\sqrt{L})$ and it satisfies
        $$
        \sup_{t,x>0} \int |\Theta_t(x,y)| d\mu(y) \leq C.
        $$
}

\pr{[Proof of Theorem \ref{Lip_Gauss}]
By the spectral theorem $T_{t} = \psi_t(\sqrt{L}) P_{\sqrt{t}}$ and
        $$
        T_t(x,y) = \int \Psi_t(x,w) P_{\sqrt{t}}(w,y) \, d\mu(w).
        $$
    Using Theorem \ref{Lip_Poiss} together with Lemma \ref{spec_mult} and Proposition \ref{Poiss_new}, for $\sqrt{t}\geq d(y,z)$, we have
        \begin{equation}
        \begin{split}
        \label{Holder}
        |T_t(x,y)-T_t(x,z)| &= \abs{\int \Psi_t(x,w)(P_{\sqrt{t}}(w,y)-P_{\sqrt{t}}(w,z)) \, d\mu(w)} \cr
        &\sleq   \eee{\frac{d(y,z)}{\sqrt{t}}}^{\a}  \int |\Psi_t(x,w)| P_{\sqrt{t}}(w,y) \, d\mu(w)\\
        &\sleq   \mu(B(y,\sqrt{t}))^{-1} \eee{\frac{d(y,z)}{\sqrt{t}}}^{\a} .
        \end{split}
        \end{equation}

We claim that for $d(y,z)\leq \sqrt{t}$ one has
    \eq{\lab{est_y}
    |T_t(x,y)-T_t(x,z)|\sleq \mu(B(y,\sqrt{t}))^{-1} \eee{\frac{d(y,z)}{\sqrt{t}}}^{\alpha\slash 2} \exp\eee{-\frac{cd(x,y)^2}{ t}}.
    }
To prove the claim we consider two cases.

{\bf Case 1:} $2d(y,z)\geq d(x,y)$. Recall that $\sqrt{t}\geq d(y,z)$, thus $d(x,y)\leq 2 \sqrt{t}$ and \eqref{est_y} follows directly from \eqref{Holder}.

{\bf Case 2:} $2d(y,z)\leq d(x,y)$. In this case $d(x,y)\simeq d(x,z)$, so by \eqref{gauss12} we obtain
    \begin{equation}
    \begin{split}\label{est123}
    |T_t(x,y) -T_t(x,z)| &\sleq  \frac{\exp\(-\frac{c_2 d(x,y)^2}{t}\)}{\mu(B(y,\sqrt{t}))}+\frac{\exp\(-\frac{c_2 d(x,z)^2}{t}\)}{\mu(B(z,\sqrt{t}))}\\
    &\sleq \mu(B(y, \sqrt{t}))^{-1} \exp\eee{-\frac{cd(x,y)^2}{t}},
    \end{split}
    \end{equation}
    where in the last inequality we have have utilized that $\mu(B(z,\sqrt{t}))\simeq \mu(B(y,\sqrt{t}))$, since $d(y,z)<\sqrt{t}$.
By taking the geometric mean of \eqref{Holder} and \eqref{est123} we obtain \eqref{est_y}. To finish the proof observe that \eqref{est_y} implies \eqref{gauss13}. Indeed,
    \alx{
    \mu(B(x,\sqrt{t}))&= \mu(B(y,\sqrt{t})) \frac{\mu(B(x,\sqrt{t}))}{\mu(B(y,\sqrt{t}))} \leq \mu(B(y,\sqrt{t})) \frac{\mu(B(y,\sqrt{t}+d(x,y)))}{\mu(B(y,\sqrt{t}))}\\
    &\leq \mu(B(y,\sqrt{t})) \eee{1+\frac{d(x,y)}{\sqrt{t}}}^q
    }
and using the exponent factor we can replace $\mu(B(y,\sqrt{t}))$ by $\mu(B(x,\sqrt{t}))$.
}

\begin{Rem}
  Let us remark that Lemma \ref{spec_mult}, which is   crucial in our  proof of  Theorem \ref{Lip_Gauss},  can be proved  by applying  functional calculus of Hebisch \cite{Hebisch_slowly}.  Thus, Theorem \ref{Lip_Gauss} can be also obtained  without using the finite speed propagation of the wave equation.
\end{Rem}

As a consequence of \eqref{gauss1212} and Theorem \ref{Lip_Gauss} we get what follows.
    \cor{cor_cont}{
    The function $T_t(x,y)$ is continuous on $(0,\8)\times X \times X$.
    }
}

As a direct consequence of Theorem \ref{Lip_Gauss} and the Doob transform (see \eqref{TK}) we get the following corollary. Notice that in Corollary \ref{Holder_Gauss} we do not assume that $T_t(x,y)$ is conservative.
\cor{Holder_Gauss}{
Assume that the semigroup $T_t$ satisfies \eqref{gauss12}. Then there are constants $\alpha,c>0$ such that
    \begin{equation*}\label{gauss1313}
    \Big|\frac{T_t(x,y)}{\w(x)\w(y)}-\frac{T_t(x,z)}{\w(x)\w(z)}\Big|\sleq \mu(B(x,\sqrt{t}))^{-1} \eee{\frac{d(y,z)}{\sqrt{t}}}^\alpha \exp\eee{-\frac{cd(x,y)^2}{ t}}
    \end{equation*}
whenever $d(y,z)\leq \sqrt{t}$.
}

\section{Measures and distances}\label{sec_new_met}

To prove Theorem \ref{mainthm2} we introduce a new quasi-metric $\wt{d}$ on $X$, which is related to $d$ and $\mu$. To this end, set
\begin{equation*}
\widetilde{d}(x,y)=\inf \mu(B),
\end{equation*}
where the infimum is taken over all closed balls  $B$ containing $x$ and $y$ (see, e.g. \cite{CoifmanWeiss_BullAMS}, \cite{Macias_Segovia_Advances_Lipschitz}). Denote
$$\wt{B}(x,r) = \set{y\in X \ | \ \wt{d}(x,y)\leq r}.$$

In the lemma below we state some properties of $\wt{d}$, which are known among specialists, and which  we shall need latter on. Since their proofs are very short and it is difficult for us to indicate one reference which contains  all of them, we provide the details for the convenience of the reader.

\lem{LemmaDistances}{The function $\wt d$ has the following properties:
\begin{itemize}
\item[(a)] there exists $C_b$ such that for $x,y\in X$ we have
\begin{equation}\label{ineqq}
C_b^{-1}\mu(B(x,d(x,y))\leq \widetilde d(x,y)\leq \mu(B(x,d(x,y)).
\end{equation}
\item[(b)]
$\,\widetilde{d}$ is a quasi-metric, namely there exists $A_1$ such that
    $$\widetilde d(x,y)\leq A_1\(\widetilde d(x,z)+\widetilde d(z,y)\).$$
\end{itemize}
Moreover, if the measure \,$\mu$ has no atoms and $\mu(X)=\8$, then:
\begin{itemize}
\item[(c)] the measure $\mu$ is regular with respect to $\wt d$, namely for $x\in X$ and $r>0$,
\begin{equation*}
\mu (\widetilde{B}(x,r))\simeq r;
\end{equation*}
\end{itemize}
\begin{itemize}
\item[(d)] for $x\in X$ and $r>0$ there exists $R>0$ such that
\begin{equation*}\label{inclusion1}
\widetilde B(x,r)\subset B(x,R)  \ \text{ and } \ \mu (B(x,R))\sleq \mu(\widetilde B(x,r));
\end{equation*}
\end{itemize}
\begin{itemize}
\item[(e)]
for $x\in X$ and $R>0$ there exists $r>0$ such that
\begin{equation*}\label{inclusion2}
B(x,R)\subset \widetilde B(x,r)  \ \text{ and } \ \mu (\widetilde B(x,r))\sleq \mu( B(x,R)).
\end{equation*}
\end{itemize}
}

\begin{proof}
(a) Set  $R=d(x,y)$. Clearly,  $\widetilde{d} (x,y)\leq \mu (B(x,R))$, as $x$ and $y$ belong to $B(x,R)$. On the other hand, if $x$ and $y$ belong to a ball $B =  B(z,r)$, then $R\leq 2r$, hence  $B(x,R) \subset B(z,3r)$ and $\mu (B(x,R))\leq \mu (B(z,3r))\simeq \mu (B(z,r))$. By taking the infimum over all balls $B$ containing both $x$ and $y$, we conclude that $\mu (B(x,R))\leq C \widetilde{d} (x,y)$ .

(b) For every  $x,y,z\in X$, we have  $d (x,y) \leq  d (x,z) +  d (z,y)$. Assume that $r=d(x,z)\geq d(z,y)$. Then $x,y\in B(z,r)$. By using (a), we deduce that
\begin{equation*}
\widetilde{d} (x,y)\leq \mu (B(z,r)) \simeq\widetilde{d} (z,x)
\leq \widetilde{d} (x,z) + \widetilde{d} (z,y)\,.
\end{equation*}

(c) Given  $x\in X$, by our additional assumptions, the function $(0,\infty) \ni r\mapsto  \mu (B(x,r))$ is increasing and
\begin{equation*}\begin{cases}
\mu (B(x,r))\searrow 0 &\text{as \;}r \searrow0\,,\\
\mu (B(x,r))\nearrow+\infty &\text{as \;}r \nearrow+\infty\,.
\end{cases}\end{equation*}
Let  $x\in X$ and  $r>0$. For every  $y\in\widetilde{B}(x,r)$, we have $\mu (B(x,d(x,y))\simeq \widetilde{d} (x,y) \le  r$ .
Hence
\begin{equation*}
R = \sup\,\{  d (x,y) |\,y\in\widetilde{B}(x,r) \}<+\infty\,.
\end{equation*}
Let  $y\in\widetilde{B}(x,r)$ such that $d(x,y) \ge \frac{R}{2}$.
Then
\begin{equation}\label{inclusion3}
\widetilde{B}(x,r)
\subset  B(x,R) \subset  B(x,2d(x,y)).
\end{equation}
Hence
\begin{equation}\label{ineq1}
\mu (\widetilde{B}(x,r))\leq \mu(B(x,R))\leq \mu (B(x,2d(x,y)))\simeq \mu (B(x,d(x,y)))\simeq \widetilde{d} (x,y)\leq r.
\end{equation}
On the other hand,
\begin{equation*}
T = \inf\{  t>0\,|\,\mu (B(x,t))\ge r \}>0\,.
\end{equation*}
As  $\mu (B(x,T/2))<r$, we have $\widetilde{d} (x,y)<r$, for every  $y\in B(x,T/2)$, hence $B(x,T/2)\subset \widetilde{B}(x,r)$.
Consequently,
\begin{equation*}\label{ineqlower}
r\le \mu (B(x,2T))\simeq \mu (B(x,T/2))
\le \mu (\widetilde{B}(x,r)) ,
\end{equation*}
which together with \eqref{ineq1} completes the proof of  (c).

(d) is a simple consequence of \eqref{inclusion3}, \eqref{ineq1}, and (c).

(e) Set $r=\mu(B(x,R))$. If $y\in B(x,R)$ then $\widetilde d(x,y)\leq r$ and, consequently, $B(x,R)\subset \widetilde B(x, r)$. Clearly, by (c),
$\mu (\widetilde B(x,r)) \simeq r = \mu(B(x,R))$.
\end{proof}

Let us recall that in Theorems \ref{mainthm1} and \ref{mainthm2} we assume that $\Phi_x$ is a bijection on $(0,\8)$. This, obviously implies that $\mu(X)=\8$ and that $\mu$ is non-atomic. As a consequence of (d) and (e) of Lemma \ref{LemmaDistances} we obtain the following corollary.

\cor{cor1123}{Suppose that $\mu$ has no atoms and $\mu(X)=\infty$. Then
the atomic Hardy spaces $\had1$ and $\hadd$ coincide and the corresponding atomic norms are equivalent.
}

We finish this section by Lemma \ref{lem1235}, which is used  latter on. Define $A_2 := C_b (C_d 2^q)^{3}$, where $C_d$, $q$, and $C_b$ are as in \eqref{doubling} and \eqref{ineqq}.

\lem{lem1235}{
Suppose that we have a space of homogeneous type $(X,d,\mu)$ such that the function $\Phi_x$ defined in \eqref{Psi} is a bijection on $(0,\8)$. Assume that $x\in X$, $r,t>0$ are related by $r=\mu(B(x,\sqrt{t}))$ and satisfy: $\sqrt{t}\leq d(y,z)$, $A_2\wt{d}(y,z)< r$.  Then
$$ \sqrt{t} \leq d(x,y)\leq  2d(x,z). $$
}

\begin{proof}
  Suppose, towards a contradiction, that $d(x,y) <\sqrt{t}$. From \eqref{ineqq} we get
  \alx{
  r&=\mu(B(x,\sqrt{t}))\leq \mu (B(y,2\sqrt{t})) \leq C_d 2^q \mu (B(y,\sqrt{t}))\\
  &\leq C_d 2^q \mu(B(y,d(y,z)))\leq C_d 2^q C_b \widetilde d(y,z)<r,
  }
  so the first inequality is proved.

  Similarly, assume $d(x,z)< d(x,y)/2$. Then $d(x,y)/2\leq d(y,z)\leq 2d(x,y)$. Thus, using \eqref{ineqq},
  \begin{equation*}\begin{split}
  \widetilde d(y,z)& \geq C_b^{-1}\mu(B(y,d(y,z))\geq C_b^{-1}(C_d 2^q)^{-1} \mu(B(y,  d(x,y))\\
  &\geq C_b^{-1}(C_d 2^q)^{-2}\mu(B(x, d(x,y)))\geq C_b^{-1}(C_d 2^q)^{-3}\mu (B(x,\sqrt{t}))=A_2^{-1} r
  \end{split}\end{equation*}
  and we come to a contradiction.
\end{proof}

\newcommand{\ddd}{\mathrm{d}}

\section{Proof of Theorem \ref{mainthm2}}\label{Uchiyama}
In order to prove Theorem \ref{mainthm2} we shall use a result of Uchiyama \cite{Uchiyama}, which we state below in Theorem \ref{Uchi}. Denote by $(X,\wt{d},\mu)$ the space $X$ equipped with a quasi-metric $\wt{d}$ and a non-negative measure $\mu$, where $\mu(X)=\8$. Assume moreover that
    \eq{\lab{metric}
     \mu(\wt{B}(x,r)) \simeq r,
    }
where $x\in X$, $r>0$ and $\wt{B}(x,r)\subseteq X$ is a ball in the quasi-metric $\wt{d}$. Let  $A_1$ be a constant in the quasi-triangle inequality, i.e.
    \eq{\label{quasi}
    \wt{d}(x,y)\leq A_1(\wt{d}(x,z)+\wt{d}(z,y)), \qquad x,y,z \in X.
    }
Additionally, we impose that there exist constants $\gamma_1, \gamma_2, \gamma_3>0$, $A\geq A_1$ and a continuous function $\wt{T}(r,x,y)$ of variables $x,y\in X$ and $r>0$ such that \al{\tag*{($U_1$)}\label{U1}
        &\wt T(r,x,x) \sgeq r^{-1},\\
        \tag*{($U_2$)}\label{U2}
        &0 \leq \wt T(r,x,y) \sleq r^{-1}\(1+\frac{\wt{d}(x,y)}{r}\)^{-1-\gamma_1},\\
        \tag*{($U_3$)}\label{U3}
        &\text{if } \wt{d}(y,z) < (r + \wt{d}(x,y))/(4A), \text{ then}\notag \\ \notag
        &\abs{\wt T(r,x,y) - \wt T(r,x,z)} \sleq r^{-1} \(\frac{\wt{d}(y,z)}{r}\)^{\gamma_2} \(1+\frac{\wt{d}(x,y)}{r}\)^{-1-\gamma_3},
    }
for all $x,y,z \in X$ and $r>0$.

As in \cite{Uchiyama}, we consider the maximal function
    \eqx{
    f^{(+)}(x) = \sup_{r>0} \abs{\int_X \wt T(r,x,y)f(y) \, d\mu(y)}
    }
and the Hardy space $H^1(X,\wt T) = \set{f\in L^1(X,\mu) \ \big| \ \norm{f}_{H^1(X,\wt T)}:= \norm{f^{(+)}}_{L^1(X,\mu)} <\8}$.

Recall that the atomic space $\hadd$ is defined in Section \ref{Intro}.

\thm{Uchi}{[\cite{Uchiyama}, Corollary 1']
Suppose that $(X,\wt{d},\mu, \wt T)$ satisfy \eqref{metric}, \eqref{quasi}, \ref{U1}, \ref{U2}, and \ref{U3}. Then the spaces $H^1(X,\wt T)$ and $\hadd$ coincide and
    \eqx{
    \norm{f}_{H^1(X,\wt T)} \simeq \norm{f}_{\hadd}.
    }
}

Assume that the kernel $T_t(x,y)$ satisfies \eqref{gauss12} and the semigroup $T_t$ is conservative. Recall that in Section \ref{sec_Lip} we proved H\"older-type estimate \eqref{gauss13} for $T_t(x,y)$. Define $\wt T(r,x,y)$ by
\eqx{\lab{KK}
 \wt T(r,x,y):=T_t(x,y),
 }
where $t=t(x,r)$ is such that
\eq{\label{tr}
\mu(B(x,\sqrt{t}))=r.
}
In what follows $t,r>0$ and $x\in X$ are always related by \eqref{tr}. Let us notice that from Corollary \ref{cor_cont} and by the assumption that $\Phi_x$ is a continuous bijection on $(0,\8)$ we have that $\wt{T}$ is a continuous function on $(0,\8)\times X \times X$.
    \thm{thm_Uchi}{
        Suppose that $T_t(x,y)$ satisfies upper and lower Gaussian bounds \eqref{gauss12} and the semigroup $T_t$ is conservative. Then the kernel $\wt T(r,x,y)$ satisfies \ref{U1}, \ref{U2}, and \ref{U3}.
    }

\pr{
The on-diagonal lower estimate \ref{U1}
is an immediate consequence of the lower Gaussian bound \eqref{gauss12}.

For every fixed   $\delta>0$ ,
the upper estimate
\begin{equation}\label{UpperEstimateKr}
\wt{T}(r,x,y)\lesssim r^{-1}\eee{ 1
+ \frac{\widetilde{d}( x, y )}r }^{ -1-\delta}
\end{equation}
follows from the upper estimates for $T_t(x,y)$,
more precisely by combining
\vspace{-1.5mm}
\begin{equation}\label{upperGauss}
\wt{T}(r,x,y)\sleq r^{-1}\exp\eee{-\frac{cd(x,y)^2}{t}}
\end{equation}
\vspace{-6.5mm}
with
\sp{\label{GaussianYieldsPolynomial}
\eee{ 1 +
\frac{\widetilde{d} (x,y)}r }^{ 1+\delta}
 &\leq \eee{ 1 +
\frac{\mu (  B
(x ,  d(   x , y) ))}
{\mu(  B
(x ,\sqrt{t } ))}
 }^{ 1+\delta} \\
 &\lesssim
\eee{1 + \frac{d(x ,y)}{\sqrt{t }}
 }^{q (1+\delta)}
\lesssim  \exp\eee{\frac{cd( x , y)^2}{t}} .
}
The latter estimate holds for any $\delta>0$. Thus \ref{U2} is proved with any $\gamma_1>0$. To finish the proof we need the H\"older-type estimate \ref{U3}. This is proved in Proposition~\ref{LemmaLipschitzEstimateKr} below.
}

\prop{LemmaLipschitzEstimateKr}{
Let $\a$ be a constant as in Theorem \ref{Lip_Gauss}. There exists  $A\geq A_1$ such that for $\delta>0$ we have
\begin{equation}\label{LipschitzEstimateKr}
\abs{\wt{T}(r, x, y) -  \wt{T}(r, x, z)}
\sleq r^{-1} \eee{ 1 + \frac{\widetilde{d} ( x,  {y})}r  }^{ -1- \delta}
\eee{ \frac{\widetilde{d}(y,  z)}r }^{ \frac{\alpha}{q}}
\end{equation}
if $\wt{d}(y, z)\leq  (r+\wt{d}(x,y))\slash (4A)$.
}

\begin{proof}
Set $A=\max(A_1,A_2)$, see \eqref{quasi} and Lemma \ref{lem1235}. Assuming that $\wt{d} (y, z) \leq  (r+\wt{d}(x,y))\slash (4A)$ let us begin with some observations.

Firstly, we claim that it suffices to prove \eqref{LipschitzEstimateKr} for $\wt d(y,z)<r\slash (2A)$.
Indeed, if $\wt d(y,z)>  r\slash (2A)$, then $\wt d(y,z)\leq \wt d(x,y)\slash (2A)$ and, consequently,
$$\wt d(x,y)\leq A_1 \wt d(x,z)+ A_1\wt d(z,y)\leq A_1\wt d(x,z)+\wt d(x,y)\slash 2.$$ So, $\wt d (x,y)\sleq \wt d (x,z)$ and
\eqref{LipschitzEstimateKr} follows from \eqref{UpperEstimateKr} by using the triangle inequality. From now on we assume that $\wt d(y,z)<r\slash (2A)$.

Secondly, if  $d(y,z)\leq \sqrt{t}$, then using \eqref{doubling} and \eqref{ineqq},
\begin{equation}\label{estimate1}
  \begin{split}
     \frac{d(y,z)}{\sqrt{t}} & \sleq \Big(\frac{\mu(B(y,d(y,z))}{\mu(B(y,\sqrt{t}))}\Big)^{1\slash q} \\
       &\leq  \Big(\frac{\mu(B(y,d(y,z))}{\mu(B(x,\sqrt{t}))}\Big)^{1\slash q}
       \Big(\frac{\mu(B(y,\sqrt{t}+d(x,y))}{\mu(B(y,\sqrt{t}))}\Big)^{1\slash q}\\
       & \sleq  \Big(\frac{\wt d(y,z)}{ r }\Big)^{1\slash q} \Big(1+\frac{d(x,y)}{\sqrt{t}}\Big).\\
  \end{split}
\end{equation}

Thirdly, if $d(y,z)\geq \sqrt{t}$, then using \eqref{doubling} and \eqref{ineqq},
\begin{equation}\label{ineq2}\begin{split}
1 & \sleq \frac{\wt d(y,z)}{\mu(B(y,\sqrt{t}))}
= \frac{\wt d(y,z)}{r}
\frac{\mu(B(x,\sqrt{t})) )}{\mu(B(y,\sqrt{t}))}\\
 &\sleq  \frac{\wt d(y,z)}{r}
\frac{\mu(B(y,\sqrt{t}+d(x,y)))}{\mu(B(y,\sqrt{t}))}\\
&\sleq \frac{\wt d(y,z)}{r} \Big(1+\frac{d(x,y)}{\sqrt{t}}\Big)^{q}.
\end{split}\end{equation}
Let us turn to the proof of \eqref{LipschitzEstimateKr}.

{\bf Case 1:} $d(y,z)\leq d(x,y)/2$. Then $d(x,y)\simeq d(x,z)$.
Thus, according to \eqref{gauss12} and \eqref{gauss13}  combined with \eqref{estimate1} and \eqref{ineq2} we obtain
\begin{equation*}\begin{split}
|  \wt{T}(r,{x},{y})
-  \wt{T}(r,{x},z) |
&\sleq
r^{-1}\,
\exp\eee{-\frac{cd(  {x} ,   {y} )^2}{t}}\,
\min\eee{\frac{d(  {y} , z)}{\sqrt{  t}}, 1}^\a\\
&\sleq r^{-1}\,
\exp\eee{-\frac{cd(  {x} ,   {y} )^2}{t}}\,
\Big(\frac{\wt d(  {y} , z)}{{  r\,}}\Big)^{\alpha \slash q} \Big(1+\frac{d(x,y)}{\sqrt{t}}\Big)^{\alpha}\\
&\sleq r^{-1}\Big(1+\frac{\wt d(x,y)}{r}\Big)^{-1-\delta}\Big(\frac{\wt d(  {y} , z)}{{  r\,}}\Big)^{\alpha \slash q},
\end{split}\end{equation*}
where in the last inequality we have used \eqref{GaussianYieldsPolynomial}.

{\bf Case 2:} $d(y,z) > d(x,y)/2$ and $\sqrt{t}>d(y,z)$. Then $\wt d(x,y)\sleq \mu(B(x,\sqrt{t}))=r$. Using \eqref{gauss13} and \eqref{estimate1},
\begin{equation*}\begin{split}
\abs{  \wt{T}(r,{x},{y})
-  \wt{T}(r,{x},z) }
&\sleq
r^{-1}\,
\eee{\frac{d(  {y} , z)}{\sqrt{  t\,}}}^{\alpha}\,
\sleq r^{-1}\eee{\frac{\wt d(  {y} , z)}{{  r\,}}}^{\alpha \slash q}.
\end{split}\end{equation*}

{\bf Case 3:} $d(y,z) >  d(x,y)/2$ and $\sqrt{t}<d(y,z)$. Then, $\wt{d}(y,z) \leq r/(2A)\leq r/(2A_2)$ and by Lemma \ref{lem1235} we have $2d(x,z)>d(x,y) \geq  \sqrt{t}$. Hence, using \eqref{upperGauss} and \eqref{ineq2},
\begin{equation*}\begin{split}
\abs{ \wt{T}(r,{x},{y})
-\wt{T}(r,{x},z)}
 & \sleq r^{-1}\exp\eee{-\frac{cd(x,y)^2}{t}}\\
 & \sleq r^{-1}\exp\eee{-\frac{cd(x,y)^2}{t}}\Big(\frac{\wt d(y,z)}{r} \Big)^{\alpha \slash q}\Big(1+\frac{d(x,y)}{\sqrt{t}}\Big)^{\alpha}\\
&\sleq  r^{-1} \Big(\frac{ \wt d(  {y} , z)}{r}\Big)^{\alpha \slash q}\Big(1+\frac{\wt d(x,y)}{r}\Big)^{-1-
\delta},
\end{split}\end{equation*}
where in the last inequality we have used \eqref{GaussianYieldsPolynomial}.
This finishes  the proof of Proposition  \ref{LemmaLipschitzEstimateKr}.
\end{proof}


Now, we gather all the elements of the proof of Theorem \ref{mainthm2}. Assuming \eqref{gauss12} and \eqref{conservative} we obtain H\"older-type estimate \eqref{gauss13}. Recall once more that the assumption on $\Phi_x$ implies $\mu(X)=\8$ and $\mu$ is non-atomic. Then we define a new quasi-metric $\wt{d}$. By Corollary \ref{cor1123} we get that $\had1=\hadd$. We apply Theorem \ref{Uchi} to the space $(X, \wt d, \mu)$. The assumptions of Theorem \ref{Uchi} are verified in Theorem \ref{thm_Uchi} and Proposition \ref{LemmaLipschitzEstimateKr}. In this way we get
$$\norm{f}_{H^1(X,\wt T)} \simeq \norm{f}_{\hadd}.$$
Using once again the assumption on $\Phi_x$ and the definition of $\wt{T}$ we easily observe that
$$\norm{f}_{H^1(X,\wt T)} = \norm{f}_{H^1(L)},$$
which finishes the proof of Theorem \ref{mainthm2}.

Let us recall that  Theorem \ref{mainthm1} follows from Theorem \ref{mainthm2}. This is elaborated at the end of  Section \ref{Doob}.

\rem{remark_ttt}{
Under the assumptions of Theorem \ref{mainthm1} one can prove, by the same methods, that the Hardy space
    \eqx{
    H^1_{\sqrt{L}} = \set{f\in L^1(X,\mu) \ \big| \ \norm{f}_{H^1_{\sqrt{L}}} := \norm{\sup_{t>0} \abs{P_tf(x)}}_{L^1(X,\mu)} <\8}
    }
coincides with $H^1_{at}(d, \w \, \mu)$ and the corresponding norms are equivalent. To this end, one uses: Proposition \ref{Poiss_new}, Doob's transform, Theorem \ref{Lip_Poiss}, and Theorem \ref{Uchi} applied to the kernel $\wt{P}(r,x,y) = P_t(x,y)$, where $t=t(x,r)$ is defined by the relation $\mu(B(x,t)) = r$.
}

\section{Examples}\label{sec_ex}

In this section we give  examples of self-adjoint semigroups   with the two-sided Gaussian bounds.

\subsection{Laplace-Beltrami operators} Let $(X,g)$ be a complete Riemannian manifold with the Riemannian distance $d(x,y)$ and the Riemannian maesure $\mu$ satisfying doubling property for balls and the Poincar\'e inequality
$$ \int_{B(x,r)} |f-f_B|^2\, d\mu \leq C r^2\int_{B(x,2r)} |\nabla f|^2\, d\mu,$$
where $\nabla$ denotes the gradient on $X$. It is well known that the kernel of the heat semigroup generated by the Laplace-Beltrami operator satisfies the two-sided Gaussian bounds \eqref{gauss12} and the H\"older estimates \eqref{gauss13}. For details and more information concerning the heat equation on Riemannian manifolds we refer the reader to \cite{Saloff-Coste} and references therein.

\subsection{Schr\"odinger operators}\label{ssec71}
On $X=\RR^d$ with the Euclidean metric and the Lebesgue measure we consider the Schr\"odinger operator
    $$L=-\Delta+V,$$
where $\Delta$ is the standard Laplacian and $V$ is a locally integrable function.

It is well-known (see, e.g. \cite{Semenov}) that for $V\geq 0$, $d\geq 3$, the semigroup $T_t=e^{-tL}$ admits kernels $T_t(x,y)$ with the upper and lower Gaussian bounds \eqref{gauss12} if and only if $V$ is a Green bounded potential, that is,
    \eq{\lab{Kato1}
    \sup_{x\in \Rd} \int_\Rd |x-y|^{2-d} V(y)\, dy <\8.
    }
Hardy spaces associated with Schr\"odinger operators on $\mathbb R^d$ satisfying \eqref{Kato1} were studied in \cite{DZ_JFAA}. Actually, as we have already mentioned, this work motivated us to study the problem of $H^1$ spaces with the Gaussian bounds in the generality as in Theorems \ref{mainthm1} and \ref{mainthm2}.

Our second example concerns Schr\"odinger operators $L=-\Delta+V$ with non-positive potentials. For $d\geq 3$ fix $p_1,p_2>1$ satisfying  $p_1<d/2<p_2$. Then there is a constant $c(p_1,p_2,d)>0$ such that if $V\leq 0$ and
$$\| V\|_{L^{p_1}(\Rd)}+\|V\|_{L^{p_2}(\Rd)}\leq c(p_1,p_2,d),$$
then the integral kernel $T_t(x,y)$ of the semigroup $T_t=e^{-tL}$ exists and satisfies two-sided  Gaussian bounds
   \eq{\lab{laba}
    T_t(x,y) \simeq t^{-d/2} \exp\eee{-\frac{|x-y|^2}{4t}}.
    }
The result can be found in Zhang \cite{Zhang_Bull_London_2003}. A~slightly different proof of \eqref{laba}, based on bridges of the Gauss-Weierstrass semigroup, can be obtained by using Lemma 1.2 together with Proposition~2.2 of \cite{Bogdan_JD_Szczypkowski}.

\subsection{Bessel-Schr\"odinger operator}

Let $\a>0$  and consider $X=(0,\8)$ and $d\mu(x)=x^\a \, dx$. Notice that the space $(X,\mu)$ with the Euclidean metric $d_e(x,y) = |x-y|$ is a space of homogeneous type. We consider the classical Bessel operator
    \eqx{
    \mathcal Bf(x) = -f''(x) - \frac{\a}{x}f'(x),
    }
which is self-adjoint positive on $L^2(X,\mu)$, and the associated semigroup of linear operators  $S_t=e^{-t\mathcal B}$. It is well-known that $S_t$ is conservative and has the integral kernel
   \eq{\label{bessel_kernel}
    S_t(x,y) = (2t)^{-1} \exp\expr{-\frac{x^2+y^2}{4t}}I_{(\a-1)/2} \expr{\frac{xy}{2t}} \expr{xy}^{-(\a-1)/2},
    }
see e.g., \cite[Chapter 6]{Borodin-Salminen}. Here $I_{(\a-1)/2}$ denotes the  modified Bessel function of the first kind, see e.g. \cite{Watson}. The kernel $S_t(x,y)$ satisfies two-sided Gaussian estimates
    \eq{\lab{gauss12S}
    \frac{C_0^{-1}}{\mu\eee{B(x,\sqrt{t})}} \exp\eee{-\frac{c_1|x-y|^2}{t}} \leq S_t(x,y)\leq \frac{C_0}{\mu\eee{B(x,\sqrt{t})}} \exp\eee{-\frac{c_2|x-y|^2}{t}}
    }
For a short proof of \eqref{gauss12S} see \cite[proof of Lemma 4.2]{DPW_JFAA}. Therefore, using Theorem \ref{mainthm2} we obtain atomic characterization of $H^1(\mathcal B)$ that was previously proved in \cite{BDT_d'Analyse}.

In this subsection we consider perturbations of $\mathcal B$ of the form
    $$
    L = \mathcal B + V,
    $$
where a potential $V$ is non-negative and locally integrable. More precisely, on $L^2(X,\mu)$ we define the quadratic form
    \eqx{
    Q(f,g) = \int_X f'(x) \overline{g'(x)}\, d\mu(x) + \int_X V(x) f(x)\overline{g(x)} \, d\mu(x),
    }
with the domain
    \eqx{
    \mathrm{Dom}\(Q\) = \set{f\in L^2(X,\mu) \ \big| \ f', \sqrt{V}f \in L^2(X,\mu)}.
    }
The form $Q$ is positive and closed. Thus, it corresponds to the unique self-adjoint operator $L$ on $L^2(X,\mu)$ with the domain
    \eqx{
    \mathrm{Dom}\(L\) = \set{f\in \rm{Dom}(Q) \ \big| \ \Big(\exists h\in L^2(X,\mu)\Big) \ \Big(\forall g\in \rm{Dom}(Q)\Big) \  Q(f,g) = \int h\,  \overline{g}\,  d\mu }.
    }
By definition, $L f = h$ when $f,h$ are related as above.

Let $T_t = \exp\eee{-tL}$ be the semigroup generated by  $-L$. The Feynman-Kac formula states that
    \eq{\lab{FK}
    T_t f(x) = E^x\eee{\exp\eee{-\int_0^t V(b_s)ds} f(b_t)},
    }
where $b_s$ is the Bessel process on $(0,\8)$ associated with $S_t$. Using \eqref{FK} one gets that the semigroup $T_t$ has the form \eqref{int_rep}, where
    \eq{
    \lab{FK2}
    0 \leq T_t(x,y) \leq S_t(x,y).
    }
Therefore the upper Gaussian estimates for $T_t(x,y)$ follows simply from \eqref{gauss12S} for any locally integrable $V\geq 0$.
On the other hand, the relation between $S_t(x,y)$ and $T_t\xy$ is given by the perturbation formula
    \eq{\label{pert}
    S_t\xy = T_t\xy + \int_0^t \int_X S_{t-s}(x,z) V(z) T_s(z,y)\, d\mu(z)\, ds
    }

From now on we consider only $\a>1$. We are interested in proving the lower Gaussian estimates, but this can be done only for some potentials $V$. For other potentials Hardy spaces may have a local character, (see e.g., \cite{Kania1}). Let
    \eqx{
    \Gamma(x,y) = \int_0^\8 S_t\xy\, dt.
    }
be the formal kernel of $\mathcal B^{-1}$. In addition to $V\in L^1_{loc}(X)$ and $V\geq 0$, we shall need one more assumption, a version of the global Kato condition, cf. \eqref{Kato1},
    \eq{\lab{Kato2}
    \norm{V}_{Kato} := \sup_{x\in X} \int_X \Gamma(x,y) V(y)\, d\mu(y) <\8.
    }
Formally, we can rephrase this as $\mathcal B^{-1}V \in L^\8(X,\mu)$. Let us point out that
 $${\Gamma}(x,y) \simeq (x+y)^{-\a+1}.$$
  This can be easily obtained from \eqref{bessel_kernel} and well-known asymptotics for the Bessel function $I_{(\alpha-1)\slash 2}$, see also \cite[Section 2]{NS_potential_ops}.

In Lemmas \ref{lem1} and \ref{lem2} we prove that under assumption \eqref{Kato2} the lower Gaussian estimates \eqref{gauss12} hold for $T_t(x,y)$. Their proofs take examples from the Euclidean setting.

\lem{lem1}{
Suppose that $\norm{V}_{Kato}<\8$. If $|x-y|\leq\sqrt{t}$, then
    \eqx{\lab{lower_local}
    T_t(x,y) \geq C_l^{-1} \mu(B(x,\sqrt{t}))^{-1}.
    }
}

\pr{
First we shall prove Lemma \ref{lem1} with an additional assumption that $\norm{V}_{Kato} \leq \e$ for a fixed small $\e>0$. By \eqref{pert} and \eqref{FK2},
    \spx{
    S_t(x,y) - T_t(x,y) =& \int_0^t \int_X S_{t-s}(x,z)V(z)T_s(z,y)d\mu(z)\, ds = \int_0^{t/2}... \ ds + \int_{t/2}^t... \ ds \cr
    \sleq& \mu(B(x,\sqrt{t}))^{-1} \int_0^{t/2}\int_X V(z) T_s(z,y)\, d\mu(z)\, ds\cr
     &+ \mu(B(y,\sqrt{t}))^{-1} \int_{t/2}^{t}\int_X S_{t-s}(x,z)V(z)\, d\mu(z)\, ds\cr
     \sleq & \mu(B(x,\sqrt{t}))^{-1} \norm{V}_{Kato}.
    }
By choosing proper $\e>0$ we deduce the thesis from the lower estimate \eqref{gauss12S} for $S_t(x,y)$.

Now assume that the norm $\norm{V}_{Kato}<\8$ is arbitrary. Fix $q>1$, such that $\norm{V}_{Kato} =q \e$. Set $V_q(x) = V(x)/q$, $\norm{V_q}_{Kato} = \e$, and let $T_t^q$ be the semigroup related to $L_q=\mathcal B+V_q$.  Using  \eqref{FK} and  H\"older's inequality,
    \sp{\label{split1}
    T_t\(\frac{\chi_{B(y,r)}(\cdot)}{\mu(B(y,r))}\)(x) &= E^x\expr{\expr{\exp\expr{-\int_0^t \frac{V(b_s)}{q}\,ds}}^q \frac{\chi_{B(y,r)}(b_t)}{\mu(B(y,r))} } \cr
        &\geq \eeee{E^x\expr{\exp\expr{-\int_0^t \frac{V(b_s)}{q}\, ds}\frac{\chi_{B(y,r)}(b_t)}{\mu(B(y,r))}}}^q         \eeee{E^x\expr{\frac{\chi_{B(y,r)}(b_t)}{\mu(B(y,r))}}}^{-\frac{q}{q'}}\cr
        & = \eeee{T_t^q\expr{\frac{\chi_{B(y,r)}(\cdot)}{\mu(B(y,r))}}(x)}^q \cdot \eeee{S_t\expr{\frac{\chi_{B(y,r)}(\cdot)}{\mu(B(y,r))}}(x)}^{-q/q'}
    }
Let us notice that
    \eqx{\lab{ae_conv}
    T_t \xy = \lim_{r\to 0} \mu\(B(y,r)\)^{-1} \int_{B(y,r)} T_t (x,z) \, d\mu(z)
    }
for a.e. $\xy$. By letting $r\to 0$ in \eqref{split1}, using \eqref{gauss12S} and the first part of the proof, for a.e. $(x,y)$ we obtain that
    \eqx{
    T_t(x,y) \geq T_t^q(x,y)^q S_t(x,y)^{-q/q'} \sgeq \mu(B(x,\sqrt{t}))^{-1}.
    }
}

\lem{lem2}{
Suppose that $\norm{V}_{Kato}<\8$. Then
    \eqx{\lab{low_gauss}
    T_t (x,y) \sgeq \mu(B(x,\sqrt{t}))^{-1} \exp\expr{-\frac{c|x-y|^2}{t}}.
    }
}

\pr{
Assume that $|x-y|^2/t \geq 1$ and set $m = \lceil 4|x-y|^2/t \rceil \geq 4$. For $i=0,... , m$, let $x_i = x + i(y-x)/m $, so that $x_0 =x$, $x_m = y$, and $|x_{i+1}-x_i| = |x-y|/m $. Denote $B_i = B(x_i, \sqrt{t}/(4\sqrt{m}))$ and observe that
\eqx{
    |y_{i+1} - y_i| \leq |y_i-x_i| + |x_i-x_{i+1}|+|x_{i+1}-y_{i+1}| \leq \frac{\sqrt{t}}{4\sqrt{m}}+\frac{\sqrt{t}}{2\sqrt{m}}+\frac{\sqrt{t}}{4\sqrt{m}}=\frac{\sqrt{t}}{\sqrt{m}}
}
for $y_i \in B_i$ and $y_{i+1} \in B_{i+1}$. Now we use the semigroup property, Lemma \ref{lem1}, and the doubling property of $\mu$ to obtain
    \spx{
    T_t \xy &= \int_X ...\int_X T_{t/m}(x,y_1) T_{t/m}(y_1,y_2)... T_{t/m}(y_{m-1},y)\, d\mu(y_1) \, ... \,d\mu(y_{m-1})  \\
    &\geq c_1^{m-1} \int_{B_1} ... \int_{B_{m-1}}\mu(B(x,\sqrt{t/m}))^{-1}... \mu(B(y_{m-1},\sqrt{t/m}))^{-1}  \, d\mu(y_1) \, ... \,d\mu(y_{m-1})  \\
    &\geq c_2^{m-1} \mu(B(x,\sqrt{t/m}))^{-1} \frac{\mu(B_1)...\mu(B_{m-1})}{\mu(B(x_1,\sqrt{t/m}))...\mu(B(x_{m-1},\sqrt{t/m}))}\\
    &\geq c_3^{m} \mu(B(x,\sqrt{t}))^{-1} = \mu(B(x,\sqrt{t}))^{-1} e^{-m\ln c_3^{-1}}\sgeq \mu(B(x,\sqrt{t}))^{-1} e^{-\frac{c|x-y|^2}{t}}.
    }
Notice that $c_1,c_2,c_3$, and $c$ in this estimate depend only on the constant $C_l$ from Lemma \ref{lem1} and the doubling constant of $\mu$.
}

Obviously, the space $(X,d_e,\mu)$ satisfies the assumptions of Theorem \ref{mainthm1}. Since we have the two-sided Gaussian estimates for $T_t$ (see Lemma \ref{lem1}, \eqref{FK2} and \eqref{gauss12S}) we obtain the following corollary.

\cor{cor_BS}{
Suppose that $\a>1$ and $V\geq 0$ satisfies \eqref{Kato2}. Then there exists $\w$ such that $0<C^{-1} \leq \w(x) \leq C$, $T_t \w = \w$ for $t>0$ and $H^1(L)=H^1_{at}(d_e, \w\, \mu)$. Moreover,
    \eqx{
     \norm{f}_{H^1(L)} \simeq \norm{f}_{H^1_{at}(d_e, \w\,\mu)}.
    }
}



{\bf Acknowledgments:} We want to thank Pascal Aucher, Fr{\'e}d{\'e}ric Bernicot, Li Chen, and Grzegorz Plebanek for their remarks. We are greatly indebt to Jacek Zienkiewicz for conversations concerning Schr\"odinger operators.


\bibliographystyle{amsplain}        
\bibliography{bib10}{}

\def\cprime{$'$}
\providecommand{\bysame}{\leavevmode\hbox to3em{\hrulefill}\thinspace}
\providecommand{\MR}{\relax\ifhmode\unskip\space\fi MR }
\providecommand{\MRhref}[2]{%
  \href{http://www.ams.org/mathscinet-getitem?mr=#1}{#2}
}
\providecommand{\href}[2]{#2}
\begin{thebibliography}{10}

\bibitem{Auscher_unpublished}
P.~Auscher, X.T. Duong, and A.~McIntosh, \emph{Boundedness of banach space
  valued singular integral operators and hardy spaces}, Unpublished preprint,
  2005.

\bibitem{Bernicot_Coulhon_Frey}
F.~Bernicot, T.~Coulhon, and D.~Frey, \emph{Gaussian heat kernel bounds through
  elliptic {M}oser iteration}, J. Math. Pures Appl. (9) \textbf{106} (2016),
  no.~6, 995--1037.

\bibitem{BDT_d'Analyse}
J.J. Betancor, J.~Dziuba{\'n}ski, and J.~L. Torrea, \emph{On {H}ardy spaces
  associated with {B}essel operators}, J. Anal. Math. \textbf{107} (2009),
  195--219.

\bibitem{Bogdan_JD_Szczypkowski}
K.~{Bogdan}, J.~{Dziuba{\'n}ski}, and K.~{Szczypkowski}, \emph{{Sharp Gaussian
  estimates for Schr\"odinger heat kernels: $L^p$ integrability conditions}},
  ArXiv e-prints (2015).

\bibitem{Borodin-Salminen}
A.N. Borodin and P.~Salminen, \emph{Handbook of {B}rownian motion---facts and
  formulae}, second ed., Probability and its Applications, Birkh\"auser Verlag,
  Basel, 2002.

\bibitem{Burkholder_Gundy_Silverstein}
D.~L. Burkholder, R.~F. Gundy, and M.~L. Silverstein, \emph{A maximal function
  characterization of the class {$H^{p}$}}, Trans. Amer. Math. Soc.
  \textbf{157} (1971), 137--153.

\bibitem{Coifman_Studia}
R.R. Coifman, \emph{A real variable characterization of {$H^{p}$}}, Studia
  Math. \textbf{51} (1974), 269--274.

\bibitem{CoifmanWeiss_BullAMS}
R.R. Coifman and G.~Weiss, \emph{Extensions of {H}ardy spaces and their use in
  analysis}, Bull. Amer. Math. Soc. \textbf{83} (1977), no.~4, 569--645.

\bibitem{Coulhon_Sikora}
T.~Coulhon and A.~Sikora, \emph{Gaussian heat kernel upper bounds via the
  {P}hragm\'en-{L}indel\"of theorem}, Proc. Lond. Math. Soc. (3) \textbf{96}
  (2008), no.~2, 507--544.

\bibitem{Davies}
E.B. Davies, \emph{Heat kernels and spectral theory}, Cambridge Tracts in
  Mathematics, vol.~92, Cambridge University Press, Cambridge, 1990.

\bibitem{DP_Argentina}
J.~Dziuba{\'n}ski and M.~Preisner, \emph{Remarks on spectral multiplier
  theorems on {H}ardy spaces associated with semigroups of operators}, Rev. Un.
  Mat. Argentina \textbf{50} (2009), no.~2, 201--215.

\bibitem{DPW_JFAA}
J.~Dziuba{\'n}ski, M.~Preisner, and B.~Wr{\'o}bel, \emph{Multivariate
  {H}\"ormander-type multiplier theorem for the {H}ankel transform}, J. Fourier
  Anal. Appl. \textbf{19} (2013), no.~2, 417--437.

\bibitem{DZ_JFAA}
J.~Dziuba{\'n}ski and J.~Zienkiewicz, \emph{On isomorphisms of {H}ardy spaces
  associated with {S}chr\"odinger operators}, J. Fourier Anal. Appl.
  \textbf{19} (2013), no.~3, 447--456.

\bibitem{DZ_Potential_2014}
\bysame, \emph{A characterization of {H}ardy spaces associated with certain
  {S}chr\"odinger operators}, Potential Anal. \textbf{41} (2014), no.~3,
  917--930.

\bibitem{Fefferman_Stein}
C.~Fefferman and E.M. Stein, \emph{{$H^{p}$} spaces of several variables}, Acta
  Math. \textbf{129} (1972), no.~3-4, 137--193.

\bibitem{Grigor'yan_HeatKernels}
A.~Grigor'yan, \emph{Heat kernel and analysis on manifolds}, AMS/IP Studies in
  Advanced Mathematics, vol.~47, American Mathematical Society, Providence, RI;
  International Press, Boston, MA, 2009.

\bibitem{Gyrya_Saloff-Coste}
P.~Gyrya and L.~Saloff-Coste, \emph{Neumann and {D}irichlet heat kernels in
  inner uniform domains}, Ast\'erisque (2011), no.~336, viii+144.

\bibitem{Hebisch_slowly}
W.~Hebisch, \emph{Functional calculus for slowly decaying kernels},
  unpublished, avaiable at http://www.math.uni.wroc.pl/$\sim \,$hebisch/.

\bibitem{Hebisch_Saloff-Coste_Grenoble}
W.~Hebisch and L.~Saloff-Coste, \emph{On the relation between elliptic and
  parabolic {H}arnack inequalities}, Ann. Inst. Fourier (Grenoble) \textbf{51}
  (2001), no.~5, 1437--1481.

\bibitem{Hofmann_Memoirs}
S.~Hofmann, G.~Lu, D.~Mitrea, M.~Mitrea, and L.~Yan, \emph{Hardy spaces
  associated to non-negative self-adjoint operators satisfying
  {D}avies-{G}affney estimates}, Mem. Amer. Math. Soc. \textbf{214} (2011),
  no.~1007, vi+78.

\bibitem{Jerison_Kenig}
D.S. Jerison and C.E. Kenig, \emph{Boundary value problems on {L}ipschitz
  domains}, Studies in partial differential equations, MAA Stud. Math.,
  vol.~23, Math. Assoc. America, Washington, DC, 1982, pp.~1--68.

\bibitem{Kania1}
E.~{Kania} and M.~{Preisner}, \emph{{Hardy spaces for Bessel-Schr{\"o}dinger
  operators}}, ArXiv e-prints (2016).

\bibitem{Latter_Studia}
R.H. Latter, \emph{A characterization of {$H^{p}({\bf R}^{n})$} in terms of
  atoms}, Studia Math. \textbf{62} (1978), no.~1, 93--101.

\bibitem{Macias_Segovia_Advances}
R.~A. Mac{\'{\i}}as and C.~Segovia, \emph{A decomposition into atoms of
  distributions on spaces of homogeneous type}, Adv. in Math. \textbf{33}
  (1979), no.~3, 271--309.

\bibitem{Macias_Segovia_Advances_Lipschitz}
\bysame, \emph{Lipschitz functions on spaces of homogeneous type}, Adv. in
  Math. \textbf{33} (1979), no.~3, 257--270.

\bibitem{NS_potential_ops}
A.~Nowak and K.~Stempak, \emph{Potential operators associated with {H}ankel and
  {H}ankel-{D}unkl transforms}, J.~Anal. Math. \textbf{131} (2017), 277--321.

\bibitem{Ouhabaz}
E.M. Ouhabaz, \emph{Analysis of heat equations on domains}, London Mathematical
  Society Monographs Series, vol.~31, Princeton University Press, Princeton,
  NJ, 2005.

\bibitem{Saloff-Coste}
L.~Saloff-Coste, \emph{Aspects of {S}obolev-type inequalities}, London
  Mathematical Society Lecture Note Series, vol. 289, Cambridge University
  Press, Cambridge, 2002.

\bibitem{Semenov}
Yu.A. Semenov, \emph{Stability of {$L^p$}-spectrum of generalized
  {S}chr\"odinger operators and equivalence of {G}reen's functions}, Internat.
  Math. Res. Notices (1997), no.~12, 573--593.

\bibitem{Song_Yan_homogeneous}
L.~{Song} and L.~{Yan}, \emph{{Maximal function characterizations for Hardy
  spaces associated to nonnegative self-adjoint operators on spaces of
  homogeneous type}}, ArXiv e-prints (2016).

\bibitem{Stein}
E.M. Stein, \emph{Harmonic analysis: real-variable methods, orthogonality, and
  oscillatory integrals}, Princeton Mathematical Series, vol.~43, Princeton
  University Press, Princeton, NJ, 1993, With the assistance of Timothy S.
  Murphy, Monographs in Harmonic Analysis, III.

\bibitem{Stein_Weiss}
E.M. Stein and G.~Weiss, \emph{On the theory of harmonic functions of several
  variables. {I}. {T}he theory of {$H^{p}$}-spaces}, Acta Math. \textbf{103}
  (1960), 25--62.

\bibitem{Uchiyama}
A.~Uchiyama, \emph{A maximal function characterization of {$H^{p}$} on the
  space of homogeneous type}, Trans. Amer. Math. Soc. \textbf{262} (1980),
  no.~2, 579--592.

\bibitem{Watson}
G.N. Watson, \emph{A treatise on the theory of {B}essel functions}, Cambridge
  Mathematical Library, Cambridge University Press, Cambridge, 1995, Reprint of
  the second (1944) edition.

\bibitem{Zhang_Bull_London_2003}
Qi~S. Zhang, \emph{A sharp comparison result concerning {S}chr\"odinger heat
  kernels}, Bull. London Math. Soc. \textbf{35} (2003), no.~4, 461--472.

\end{thebibliography}

\end{document}